# ACCELERATED CONVERGENCE FOR NONPARAMETRIC REGRESSION WITH COARSENED PREDICTORS

By Aurore Delaigle,[1] Peter Hall[2] and Hans-Georg Müller[3]

*University of Bristol and University of Melbourne,*
*University of Melbourne and University of California, Davis*

We consider nonparametric estimation of a regression function for a situation where precisely measured predictors are used to estimate the regression curve for coarsened, that is, less precise or contaminated predictors. Specifically, while one has available a sample $(W_1, Y_1), \ldots, (W_n, Y_n)$ of independent and identically distributed data, representing observations with precisely measured predictors, where $E(Y_i | W_i) = g(W_i)$, instead of the smooth regression function $g$, the target of interest is another smooth regression function $m$ that pertains to predictors $X_i$ that are noisy versions of the $W_i$. Our target is then the regression function $m(x) = E(Y | X = x)$, where $X$ is a contaminated version of $W$, that is, $X = W + \delta$. It is assumed that either the density of the errors is known, or replicated data are available resembling, but not necessarily the same as, the variables $X$. In either case, and under suitable conditions, we obtain $\sqrt{n}$-rates of convergence of the proposed estimator and its derivatives, and establish a functional limit theorem. Weak convergence to a Gaussian limit process implies pointwise and uniform confidence intervals and $\sqrt{n}$-consistent estimators of extrema and zeros of $m$. It is shown that these results are preserved under more general models in which $X$ is determined by an explanatory variable. Finite sample performance is investigated in simulations and illustrated by a real data example.

Received May 2005; revised February 2007.

[1]Supported by a Hellman Fellowship and by a Belgian American and Educational Foundation Post-Doctoral Fellowship, held at Department of Statistics, University of California, Davis.

[2]Supported in part by an Australian Research Council Fellowship.

[3]Supported in part by National Science Foundation Grants DMS-02-04869, DMS-03-54448 and DMS-05-05537.

## 1. Introduction.

1.1. *Motivation and models.* In this paper, we consider nonparametric estimation of a regression function in the framework of a novel errors-in-variables problem. In the classical errors-in-variables problem, the interest is to estimate a regression function $m$, where

$$Y = m(G) + \epsilon,$$

and a sample $(F_1, Y_1), \ldots, (F_n, Y_n)$ of independent and identically distributed (i.i.d.) data is available, with $F_i = G_i + \delta_i$, where $G$ and $\delta$ are independent random variables and the distribution of $\delta$ is known. References include Fan, Truong and Wang [9], Fan and Masry [7], Fan and Truong [8], Stefanski and Cook [16], Carroll, Ruppert and Stefanski [4], Carroll, Maca and Ruppert [3], Taupin [17], Devanarayan and Stefanski [5], Ioannides and Matzner-Løber [12], Linton and Whang [14] and Carroll and Hall [2].

The situation we consider here is different: we assume that an i.i.d. sample $(W_1, Y_1), \ldots, (W_n, Y_n)$ is observed, where

$$(1.1) \qquad\qquad Y_i = g(W_i) + \varepsilon_i \qquad \text{for } 1 \le i \le n,$$

with independent errors $\varepsilon_i$ with mean zero and finite variance. Instead of estimating the regression function $g(w) = \mathrm{E}(Y|W = w)$ generating the observations, the goal is to estimate the target regression function $m(x) = \mathrm{E}(Y|X = x)$, which differs from $g$, as $X$ is a contaminated (coarsened) version of $W$.

Specifically, $X \sim f_X$ and $X = W + \delta$, where $\delta \sim f_\delta$ represents a random distortion, and $W$ and $\delta$ are independent random variables. We refer to $X$ as a coarsened predictor of $Y$. In Section 1.3 we shall note that the model for $X$ can be generalized, without altering the main properties of our methods, to the situation where $X$ is a proxy for a variable $T$ related to $W$, provided we have additional data to infer the relationship between $T$ and $X$.

The motivating idea is that the sample $(W_1, Y_1), \ldots, (W_n, Y_n)$, where one has precise predictors, is hard to obtain, and therefore future values of $Y$ will be predicted from easier-to-obtain contaminated observations $X$ of $W$. This type of problem arises in situations where it is expensive or involved to measure $W$ accurately, so that, in routine applications, only the contaminated and less precise predictors $X$ are available. At the same time, a training set is available containing more precise predictors. For example, if we have a sample of repeated contaminated observations of the predictor for several individuals, the averaged observations $W_i = \bar{X}_{i\cdot}$ will provide relatively accurate measurements of the predictor.

The problem we address is how to use the information in the training sample, with its accurate measurements, to predict a future response $Y$ from a future contaminated predictor $X$. One of our central findings is that this coarsening of the predictor has the consequence of accelerating the convergence of the proposed estimator of $m$ from the usual nonparametric rate, strictly slower than $\sqrt{n}$, to a parametric $\sqrt{n}$-rate, even if the target regression function is known only to be smooth and does not follow any particular parametric model.

In the setting of (1.1), $m$ is generally not identifiable unless we know $f_\delta$. The latter assumption is commonly made in errors-in-variables problems. See, for example, Stefanski and Carroll [15] and Fan [6]. However, if we have additional data directly on $\delta$, or if the data at (1.1) are replicated, then we can identify $m(x)$ without knowing $f_\delta$. In either of these settings we might have a parametric model for $f_\delta$, or we might wish to treat inference about $f_\delta$ from a nonparametric viewpoint. In order to show that estimation of $m$ is a semiparametric problem, even if $f_\delta$ is not known and we treat it nonparametrically, we shall consider a more general, relatively "uninformative" type of replication, where we observe only

$$(1.2) \qquad U_{ij} = V_i + \delta_{ij} \qquad \text{for } 1 \le j \le n_i \text{ and } 1 \le i \le N.$$

Here, $V_1, \ldots, V_N$ are arbitrary random variables, $\delta_{11}, \ldots, \delta_{Nn_N}$ are mutually independent, the $\delta_{ij}$ are all distributed as $\delta$, and it is assumed that each $n_i \ge 2$. Our results demonstrate that it is possible to attain $\sqrt{n}$-consistency without making joint assumptions about the data at (1.1) and (1.2). In particular, it is not necessary to suppose that the $U_{ij}$ are independent of the $(W_i, \delta_i)$ or that the $V_i$ are independent of the $\delta_{ij}$. A direct application of the model in (1.2) is where $U_{ij}$ are replicated measurements of $X_i$, and $V_i = W_i$.

1.2. *Estimators.* First we express $m$ as a ratio, where each component can be estimated separately. Since $m(x) = \mathrm{E}(Y|X = x) = \mathrm{E}(g(W)|X = x)$, then

$$
\begin{aligned}
m(x) &= \frac{\int g(w) f_{X|W}(x|w) f_W(w)\, dw}{f_X(x)} \\
&= \frac{\int g(w) f_\delta(x-w) f_W(w)\, dw}{\int f_\delta(x-w) f_W(w)\, dw} = \frac{\varphi(x)}{\psi(x)},
\end{aligned}
$$

where we define $\psi(x) = \int f_\delta(x-w) f_W(w)\, dw = \mathrm{E}(f_\delta(x-W))$ and

$$\varphi(x) = \int g(w) f_\delta(x-w) f_W(w)\, dw = \mathrm{E}(g(W) f_\delta(x-W)) = \mathrm{E}(Y f_\delta(x-W)).$$

If the data $(W_i, Y_i)$ are generated by the model (1.1), and $f_\delta$ is assumed known, then the representations above motivate the estimators

$$\hat\varphi(x) = n^{-1} \sum_{i=1}^n Y_i f_\delta(x - W_i),$$

$$\hat{\psi}(x) = n^{-1} \sum_{i=1}^{n} f_\delta(x - W_i)$$

of $\varphi(x)$ and $\psi(x)$, respectively, leading to the estimators

$$(1.4) \qquad \hat{m}(x) = \frac{\sum_{i=1}^{n} Y_i f_\delta(x - W_i)}{\sum_{i=1}^{n} f_\delta(x - W_i)} = \frac{\hat{\varphi}(x)}{\hat{\psi}(x)}$$

of $m(x)$. An attractive feature of $\hat{m}$ is that it does not require a smoothing parameter.

When additional data following (1.2) are available, we propose a Fourier-inversion approach to estimating $\varphi$ and $\psi$, as follows. Assume that $\delta$ has a symmetric distribution, with positive characteristic function $f_\delta^{\mathrm{ft}}$,

$$(1.5) \qquad f_\delta^{\mathrm{ft}}(t) = \Re f_\delta^{\mathrm{ft}}(t) > 0 \qquad \text{for all real } t,$$

where the superscript ft denotes Fourier transform, and the Fourier transform of a function $f$ is given by $f^{\mathrm{ft}}(t) = \int f(x) e^{itx} \, dx$. The real part of $f^{\mathrm{ft}}$ is denoted by $\Re f^{\mathrm{ft}}$. (Our methods can be generalized to the case of asymmetric error distributions, using techniques borrowed from Li and Vuong [13].) Our estimator of $f_\delta^{\mathrm{ft}}$ is

$$(1.6) \qquad \hat{f}_\delta^{\mathrm{ft}}(t) = \left| \frac{1}{M} \sum_{j=1}^{N} \sum_{1 \le k_1 < k_2 \le n_j} \exp[it(U_{jk_1} - U_{jk_2})] \right|^{1/2},$$

where $M = \frac{1}{2} \sum_{j=1}^{N} n_j(n_j - 1)$. (Here and below, $\hat{f}^{\mathrm{ft}}$ denotes an estimator of the Fourier transform of $f$, not the Fourier transform of an estimator $\hat{f}$ of $f$.)

Writing $f_W$ for the density of $W$, estimators of the Fourier transforms, $f_W^{\mathrm{ft}}$ and $(f_W g)^{\mathrm{ft}}$, of $f_W$ and $f_W g$ are respectively given by

$$(1.7) \quad \hat{f}_W^{\mathrm{ft}}(t) = \frac{1}{n} \sum_{j=1}^{n} \exp(itW_j), \qquad \widehat{(f_W g)}^{\mathrm{ft}}(t) = \frac{1}{n} \sum_{j=1}^{n} Y_j \exp(itW_j).$$

Estimators of $\psi$ and $\varphi$ based on Fourier inversion are then obtained as

$$(1.8) \qquad \begin{aligned} \tilde{\psi}(x) &= \frac{1}{2\pi} \int_{|t| \le \tau_n} \hat{f}_W^{\mathrm{ft}}(t) \hat{f}_\delta^{\mathrm{ft}}(t) e^{-itx} \, dt, \\ \tilde{\varphi}(x) &= \frac{1}{2\pi} \int_{|t| \le \tau_n} \widehat{(f_W g)}^{\mathrm{ft}}(t) \hat{f}_\delta^{\mathrm{ft}}(t) e^{-itx} \, dt, \end{aligned}$$

where $\tau_n$ is a smoothing parameter. Our estimator of $m$ is $\tilde{m} = \Re\tilde{\varphi}/\Re\tilde{\psi}$.

1.3. *Generalizations.* The main features of our approach also apply to the more general case where $X = p(T \mid \theta) + \delta$, that is, where $W = p(T \mid \theta)$ for a r.v. $T$, and $(T, Y)$ rather than $(W, Y)$ is observed in a subset of the available data. Here $p(\cdot \mid \theta)$ is a parametric model, determined by the finite parameter $\theta$.

In this setting, we would ideally take $W_i = p(T_i \mid \theta)$. However, in most cases, we have to settle instead for $\widehat{W}_i = p(T_i \mid \hat{\theta})$, where $\hat{\theta}$ is a $\sqrt{n}$-consistent estimator of $\theta$, computed by least squares from data $(T_i', X_i')$, with the same distribution as $(T, X)$, and related by $X_i' = p(T_i' \mid \theta) + \delta_i'$, for $1 \le i \le r$, say. The most important special case is that where $p$ is linear: $p(t \mid \theta) = \theta^{(1)} + \theta^{(2)} t$, with $\theta = (\theta^{(1)}, \theta^{(2)})$ denoting a vector of length 2.

In this model, the variable $X$ typically represents a proxy for the variable $T$, where $T$ often is not available in applications, because it is too costly to measure it, for example. In some applications, however, we are able to observe $(T_i, X_i, Y_i)$ for $1 \le i \le n$ in a "training set," where $r = n$. We then propose to use the estimator $\hat{m}$, rather than a more conventional nonparametric regression based on $(X_i, Y_i)$, since it is more accurate, as we will demonstrate. In some cases the training set $(T_i', X_i')$ might be genuinely different from $(T_i, X_i)$. For example, $(T_i', X_i')$ might represent external data.

In the case where $X = p(T \mid \theta) + \delta$ the estimators $\hat{m}$ and $\tilde{m}$ differ only in that we replace $W_i$ by $\widehat{W}_i$ at each appearance. Under appropriate regularity conditions the main properties of $\hat{m}$ and $\tilde{m}$, and in particular their $\sqrt{n}$-consistency [provided $n = O(r)$], do not change. This point will be discussed in Section 2.

## 2. Asymptotic results.

2.1. *Case where $f_\delta$ is assumed known.* Here we discuss asymptotic properties of the estimator defined at (1.4). A central result is the weak convergence of a suitably scaled estimator process, with $\sqrt{n}$-scaling, to a Gaussian limit process in the location argument $x$. This result (Theorem 1 below) implies, among other matters, pointwise and uniform limits, local and simultaneous confidence bands, and convergence of estimated extrema locations.

We assume throughout Section 2 that the distribution corresponding to $f_\delta$ is absolutely continuous, and in Section 2.1 that $f_\delta$ has a bounded derivative. In Section 2.1 it is not necessary to suppose that the densities $f_W$ or $f_{W,Y}$ exist, although it is convenient to use the notation $f_W$ and $f_{W,Y}$ when introducing the quantities needed to state and derive our results. However, the differential elements $f_W(w) \, dw$ and $f_{W,Y}(w, y) \, dw \, dy$ may be interpreted as $F_W(dw)$ and $F_{W,Y}(dw, dy)$, respectively; the distributions need not be absolutely continuous.

Given an integer $\nu \geq 0$, and assuming all quantities are well defined, let $\varphi$ and $\psi$ be as in Section 1 and define

$$h(x, w) = f_\delta^{(\nu)}(x - w),$$

$$\alpha(x) = \varphi^{(\nu)}(x) = \iint y f_{W,Y}(w, y) h(x, w) \, dw \, dy,$$

$$\beta(x) = \psi^{(\nu)}(x) = \int h(x, w) f_W(w) \, dw.$$

Below, the notation $D$ denotes a compact set on which we shall estimate $\varphi$ and $\psi$. The following conditions, indexed by $\nu = 0$ or 1, will be used to prove our results. For $\nu = 1$ they can be relaxed to an assertion about the modulus of continuity for the corresponding quantity when $\nu = 0$; we impose the more stringent condition only for simplicity and brevity.

($A_{\nu,1}$) (boundedness of $f_\delta^{(\nu)}$) $\sup_{x,y \in \mathbb{R}} |h(x, y)| < \infty$;

($A_{\nu,2}$) (smoothness of $f_\delta^{(\nu)}$) $h(x, w)$ is an integrable function which is uniformly Lipschitz continuous in $x$, that is, $\sup_w |h(x_1, w) - h(x_2, w)| \leq L|x_1 - x_2|$, for a constant $L > 0$;

($A_{\nu,3}$) (boundedness of $\beta^{-1}$) $\inf_{x \in D} |\beta(x)| = c_\beta > 0$;

($A_4$) (finiteness of moments) $\int |y| f_Y(y) \, dy < \infty$ and $\int y^2 f_Y(y) \, dy < \infty$.

Note in particular that conditions ($A_{\nu,1}$) and ($A_4$) guarantee that all the quantities defined above exist, and $\alpha$ and $\beta$ satisfy $\sup_{x \in D} |\alpha(x)| < \infty$ and $\sup_{x \in D} |\beta(x)| < \infty$.

Let $\Rightarrow$ denote weak convergence in $\mathcal{C}(D)$ and define

$$\mu(x_1, x_2) = \int y f_{W,Y}(w, y) f_\delta(x_1 - w) f_\delta(x_2 - w) \, dw,$$

$$\varphi_1(x_1, x_2) = \int y^2 f_{W,Y}(w, y) f_\delta(x_1 - w) f_\delta(x_2 - w) \, dw,$$

$$\psi_1(x_1, x_2) = \int f_W(w) f_\delta(x_1 - w) f_\delta(x_2 - w) \, dw.$$

Our main result is a functional limit theorem for the proposed estimator. (All proofs are deferred to Section 5.)

THEOREM 1. *Under conditions* ($A_{\nu,1}$), ($A_{\nu,2}$) *for* $\nu = 0, 1$, ($A_{0,3}$) *and* ($A_4$), *we have that, for the process* $Z_n(x) = \sqrt{n}(\hat{m}(x) - m(x))$, $Z_n \Rightarrow Z$ *on* $\mathcal{C}(D)$, *where* $Z$ *is a Gaussian process with zero mean and covariance*

$$\text{cov}\{Z(x_1), Z(x_2)\}$$
$$= \varphi_1(x_1, x_2) / \{\psi(x_1)\psi(x_2)\} + \psi_1(x_1, x_2)\varphi(x_1)\varphi(x_2) / \{\psi^2(x_1)\psi^2(x_2)\}$$
$$- \mu(x_1, x_2)\{\varphi(x_1)\psi(x_2) + \varphi(x_2)\psi(x_1)\} / \{\psi^2(x_1)\psi^2(x_2)\},$$

*for* $x_1, x_2 \in D$.

The correlation structure for estimates at points $x_1 \neq x_2$ is seen not to vanish asymptotically, in contrast to the well-known behavior of local smoothing estimators where estimates at different points become asymptotically uncorrelated as bandwidths and windows converge to zero. Define $\hat{\mu}(x_1, x_2) = n^{-1} \sum_{i=1}^{n} Y_i f_\delta(x_1 - W_i) f_\delta(x_2 - W_i)$, $\hat{\psi}_1(x_1, x_2) = n^{-1} \sum_{i=1}^{n} f_\delta(x_1 - W_i) f_\delta(x_2 - W_i)$ and $\hat{\varphi}_1(x_1, x_2) = n^{-1} \sum_{i=1}^{n} Y_i^2 f_\delta(x_1 - W_i) f_\delta(x_2 - W_i)$. Particular consequences of Theorem 1 include the properties $\sup_{x \in D} \sqrt{n}(\hat{m}(x) - m(x)) \overset{D}{\to} \sup_{x \in D} Z(x)$ and $\sqrt{n}(\hat{m}(x) - m(x)) \overset{D}{\to} N(0, V(x))$ as $n \to \infty$, where $V(x) = \text{cov}(Z(x), Z(x))$ is estimated uniformly and $\sqrt{n}$-consistently by $\hat{V}(x) = \hat{\varphi}_1(x, x)\hat{\psi}^{-2}(x) + \hat{\varphi}^2(x)\hat{\psi}_1(x, x)\hat{\psi}^{-4}(x) - 2\hat{\varphi}(x)\hat{\mu}(x, x)\hat{\psi}^{-3}(x)$, in the sense that $\sup_{x \in D} |\hat{V}(x) - V(x)| = O_P(n^{-1/2})$. It follows that an asymptotic $(1 - \alpha)$-level confidence interval for $m(x)$ has endpoints $\hat{m}(x) \pm \hat{V}(x)^{1/2}\Phi^{-1}(1 - \alpha/2), \hat{m}(x)$, where $\Phi$ denotes the standard normal distribution function.

Semiparametric efficiency of $\hat{m}$ can be established, in regular cases where $f_\delta(x - w)$ is monotone in $x$ for $w$ in the support of $W$, by considering the following simpler problem. Suppose we observe independent and identically distributed pairs $(R_1, S_1), \ldots, (R_n, S_n)$, where $R_i \geq 0$ and $S_i = \rho(R_i) + \epsilon_i$, with $\rho$ a smooth function and $\epsilon_i$ independent of $R_i$ and distributed as $N(0, \sigma^2)$. Consider the problem of estimating $(\theta_1, \theta_2) = (E(R), E\{R\rho(R)\})$ from these data. The estimator $(\hat{\theta}_1, \hat{\theta}_2) = (n^{-1} \sum_{i=1}^{n} R_i, n^{-1} \sum_{i=1}^{n} R_i S_i)$ is asymptotically normally distributed and semiparametric efficient in this problem, and thus $\hat{\theta}_2 / \hat{\theta}_1$ is semiparametric efficient for $\theta_2 / \theta_1$. (The proof follows via Examples 3.2.1 and 3.3.4, and Propositions 3.3.1 and A.5.2, of Bickel et al. [1].) We may identify $m(x)$ and $\hat{m}(x)$ with $\theta_2 / \theta_1$ and $\hat{\theta}_2 / \hat{\theta}_1$, respectively, by taking $R_i = f_\delta(x - W_i)$ and $\rho(r) = g\{f_{\delta x}^{-1}(r)\}$, where $f_{\delta x}(w) = f_\delta(x - w)$.

Under additional regularity conditions, Theorem 1 continues to hold, although with an altered covariance structure for the limiting process $Z$, in the more general setting described in Section 1.3. There we observe $T_i$, in the setting of an unknown parameter $\theta$, rather than $W_i = p(T_i \mid \theta)$, and $W_i$ is replaced by $\widehat{W}_i = p(T_i \mid \hat{\theta})$ in the definition of $\hat{m}$. If the model $p(\cdot \mid \theta)$ is linear, then the only additional assumptions needed are two bounded derivatives of $f_\delta$, and $E(T^2) < \infty$, where $T$ denotes a generic $T_i$. See Section 5.3 for an outline proof.

2.2. *Case where $f_\delta$ is estimated from replicated data.* The conditions imposed below [see particularly (2.2)] imply that the distributions of $W$ and $\delta$ are absolutely continuous, and in particular that the respective densities $f_W$ and $f_\delta$ are square-integrable. We shall assume that

$$(2.1) \quad \max_{i \geq 1} n_i < \infty, \qquad n = o(N), \qquad n^{1/(2(\lambda + \lambda_\delta - 1))} \ll \tau_n \ll N^{1/(2(1 + \lambda_\delta))},$$

where $a_n \ll b_n$ for positive sequences $a_n$ and $b_n$ means that $a_n/b_n \to 0$ as $n \to \infty$; and that, for constants $\lambda, \lambda_\delta > 0$ satisfying

$$(2.2) \qquad\qquad \lambda > \lambda_\delta + 1 \quad \text{and} \quad \lambda_\delta > 1,$$

we have

$$(2.3) \qquad |(f_W g)^{\text{ft}}(t)| + |f_W^{\text{ft}}(t)| \leq \text{const.} |t|^{-\lambda} \qquad \text{for all } t,$$

$$f_\delta^{\text{ft}}(t) > 0 \qquad \text{for all } t, \qquad |f_\delta^{\text{ft}}(t)| \asymp |t|^{-\lambda_\delta} \qquad \text{as } t \to \infty.$$

The second part of (2.1) asks that there be an order of magnitude more values of $U_{ij}$, at (1.2), than there are pairs $(W_i, Y_i)$, at (1.1). Conditions (2.2) and (2.3) ask that $f_\delta$ be sufficiently smooth, with its Fourier transform decaying in the standard polynomial way, and that $f_W$ and $f_W g$ be sufficiently smooth relative to $f_\delta$. The second part of (2.1), and (2.2), imply that it is always possible to choose the smoothing parameter $\tau_n$ such as to satisfy the third part of (2.1).

THEOREM 2. *If the function $g$ is uniformly bounded, if the errors $\epsilon_i$ at (1.1) have zero mean and finite variance, and if (2.1)–(2.3) hold, then, uniformly in $x$,*

$$(2.4) \qquad \tilde{\psi}(x) = \hat{\psi}(x) + o_p(n^{-1/2}), \qquad \tilde{\varphi}(x) = \hat{\varphi}(x) + o_p(n^{-1/2}).$$

Let $\mathcal{I}$ denote an interval for which $\inf_{x \in \mathcal{I}} \psi(x) > 0$. Result (2.4) implies that, under the additional conditions imposed for Theorem 1, the estimator $\tilde{m} = \tilde{\varphi}/\tilde{\psi}$, which is an alternative to $\hat{m} = \hat{\varphi}/\hat{\psi}$ discussed in Section 2.1, satisfies $\tilde{m}(x) = \hat{m}(x) + o_p(n^{-1/2})$ uniformly in $x \in \mathcal{I}$. Therefore $\tilde{m}$ inherits the weak convergence and semiparametric-efficiency properties of $\tilde{m}$ on $\mathcal{I}$. Theorem 2 holds, under more restrictive assumptions, in the more general setting of Section 1.3; see Section 5.3.

**3. Simulations.** We implemented our estimator $\hat{m}(x)$ of $m(x)$ on samples of $(W, Y)$ generated from models of two types:

(1) $g(w) = [3w + 20(2\pi)^{-1/2} \exp(-200(w - 1/2)^2)]1_{[0,1]}(w)$, $W \sim U[0,1]$, $\epsilon \sim \text{N}(0, \sigma_\epsilon^2)$ and $\delta \sim \text{N}(0, \sigma_\delta^2)$ or $\delta \sim U[-1/2, 1/2]$;

(2) $Y|W = w \sim \text{Bernoulli}(g(w))$, with $g(w) = \exp(6w)/[1 + \exp(6w)]$, $W \sim U[-0.5, 0.5]$, $\delta \sim \text{N}(0, \sigma_\delta^2)$ or with $g(w) = 0.45\sin(a\pi w) + 0.5$, $a = 2$ or $4$, $W \sim U[0,1]$, $\delta \sim \text{N}(0, \sigma_\delta^2)$ or $\delta \sim U[-1/2, 1/2]$.

The last example was used by Hobert and Wand [11]. In each case, we considered several sample sizes ($n = 50, 100$ and $250$) and the parameters var($\delta$) and var($\epsilon$) were chosen such that the noise-to-signal ratios $NS_\delta = \text{var}(\delta)/\text{var}(W)$ and $NS_\epsilon = \text{var}(\epsilon)/\|g\|_\infty$ equal 10%, 25% or 50%. We considered the situation where the values of $X$ are available as well, which allowed us to compare our estimators with the $n^{-2/5}$-consistent Nadaraya–Watson estimator $\hat{m}_N$ of $m(x)$, based on observations of $(X, Y)$. In all cases, our estimators based on $(W, Y)$ performed much better than $\hat{m}_N$, which was biased and much more variable. These findings continued to hold in the setting of Section 1.3, where the sample available was $(T_i, X_i, Y_i)$, $i = 1, \ldots, n$, and the error variance was unknown and estimated by the empirical variance of the sample $X_i - \hat{W}_i$, $i = 1, \ldots, n$. More details are available from the first author's website.

The typical behavior of our estimator is illustrated in Figure 1, where we compare, for case (1) with uniform $\delta$, $NS_\delta = 0.1$, $NS_\delta = 0.25$ and $n = 250$, the results of 1000 replications of the estimators $\hat{m}$ with the correct error density $f_\delta$ and $\hat{m}$ with $f_\delta$ misspecified (here we used Gaussian error instead of the uniform error). In both cases, the estimates shown correspond to the first, fifth and ninth deciles of the ordered 1000 values of the integrated squared error $\int (\hat{m}(x) - m(x))^2 \, dx$. We see that for small $NS_\delta$, the estimator is quite robust to error misspecification, but, without any surprise, the quality deteriorates as the ratio increases. Note, however, that the results remain quite good for $NS_\delta = 0.25$.

## 4. Real data illustration.

We illustrate the proposed estimator in the setting of Section 1.3 on a real data example. The data set was collected during a South African study on heart disease and was used by Hastie, Tibshirani and Friedman [10]. The data are available at www-stat.stanford.edu/ElemStatLearn. During the study, several variables were measured on males in a heart-disease high-risk region of the Western Cape, including low density lipoprotein cholesterol (LDL) and total cholesterol (CHOL) as predictors, and coronary heart disease history (CHD) as response, coded as 0 = nonincidence of CHD, 1 = incidence of CHD. LDL is much more difficult to measure than CHOL, which motivates the use of CHOL as a proxy for LDL (Carroll, Ruppert and Stefanski [4]). After deleting several outliers, the relationship between LDL and CHOL can be reasonably well modeled as $\log(\text{CHOL}) = \theta^{(1)} + \theta^{(2)} \log(\text{LDL}) + \delta$ with $\delta$ a random variable of zero mean; see Carroll, Ruppert and Stefanski [4], who use the same model for a similar data set. Checking for outliers, we deleted the observations corresponding to the smallest (resp., two largest) value(s) of CHOL, the smallest three (resp., largest two) values of LDL, and the eight points of $(\log(\text{CHOL}), \log(\text{LDL}))$ the furthest away from the least squares line.

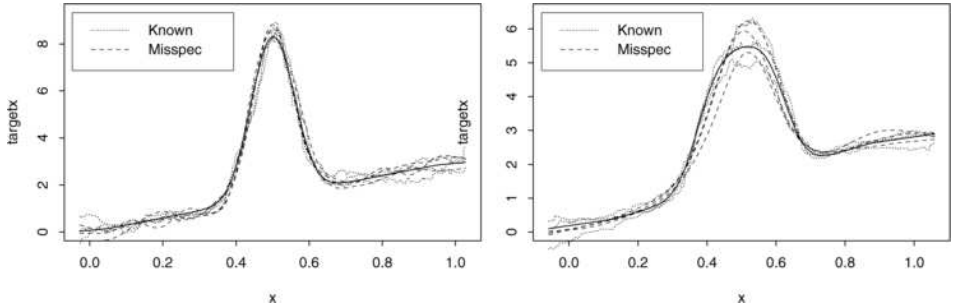

FIG. 1. *The estimator $\hat{m}$ with the error $f_\delta$ known (uniform) or misspecified (Gaussian) for case (2), with $NS_\delta = 0.1$ (left panel) or $NS_\delta = 0.25$ (right panel), with $NS_\epsilon = 0.1$ and $n = 250$. The solid curve is the target curve $m$.*

We set $Y = \text{CHD}$, $X = \log(\text{CHOL})$ and $\hat{W} = \hat{\theta}^{(1)} + \hat{\theta}^{(2)} \log(\text{LDL})$, where $\hat{\theta}^{(1)} = 4.8890$ and $\hat{\theta}^{(2)} = 0.3663$ are the least squares estimators of $\theta^{(1)}$ and $\theta^{(2)}$. Our goal is to estimate $m(x) = \text{E}(Y|X = x)$, the conditional expectation of incidence of coronary heart disease given the (transformed) total cholesterol level, using the sample of $n = 446$ observations.

We compare the proposed estimator $\hat{m}(x)$ with the Nadaraya–Watson estimator $\hat{m}_N$. The data suggest that it is reasonable to assume that the errors $\delta_i = X_i - W_i$ are normal, where the variance can be estimated from the differences $X_i - \hat{W}_i$. In Figure 2, we overlay the proposed estimator $\hat{m}$ and the Nadaraya–Watson estimator $\hat{m}_N$ calculated with an appropriate data-driven cross-validation bandwidth. The graphs suggest that the probability of coronary heart disease increases with the cholesterol level. The increase is highly nonlinear, and there are clear differences between the classical Nadaraya–Watson estimator and the proposed estimator. The Nadaraya–Watson estimator exhibits additional fluctuations, especially in the right tail, thus giving a less stable appearance.

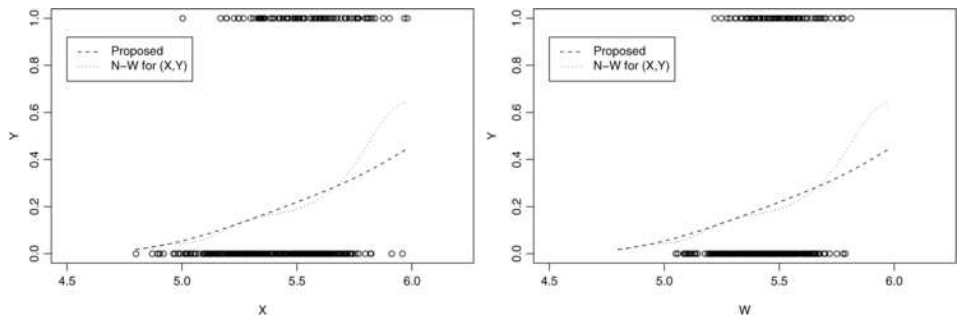

FIG. 2. *The proposed estimator $\hat{m}$ and the Nadaraya–Watson estimator $\hat{m}_N$ based on the observations of $(X, Y)$, and a scatter plot of the 446 observed values of $(X, Y)$ (left panel) or the 446 observed values of $(\hat{W}, Y)$ (right panel), for the coronary heart disease data.*

**5. Proofs.**

5.1. *Outline proof of Theorem* 1. Define the auxiliary quantities

$$\tilde{Z}_n(x) = \sqrt{n}\left[\frac{\hat{\varphi}(x) - \varphi(x)}{\psi(x)} - \frac{(\hat{\psi}(x) - \psi(x))\varphi(x)}{\psi^2(x)}\right],$$

$$\hat{\alpha}(x) = n^{-1}\sum_{i=1}^n Y_i h(x, W_i),$$

$$\alpha_1(x_1, x_2) = \iint y^2 h(x_1, w) h(x_2, w) f_{W,Y}(w, y)\, dw\, dy,$$

$$\hat{\beta}(x) = n^{-1}\sum_{i=1}^n h(x, W_i),$$

$$\beta_1(x_1, x_2) = \int h(x_1, w) h(x_2, w) f_W(w)\, dw.$$

The next two results will be useful to prove the theorem. Their proof is given at the end of this section.

LEMMA 1. *Let $\nu$ be a positive integer and $x \in D$. Under Conditions* $(A_{\nu,1})$, $(A_{\nu,2})$ *and* $(A_4)$,

$$\sqrt{n}(\hat{\alpha}(x) - \alpha(x)) \Rightarrow Z_\alpha(x), \qquad \sqrt{n}(\hat{\beta}(x) - \beta(x)) \Rightarrow Z_\beta(x),$$

*where $Z_\alpha, Z_\beta$ are Gaussian processes characterized by the moments $\mathrm{E}(Z_\alpha(x)) = \mathrm{E}(Z_\beta(x)) = 0$, and $\mathrm{cov}(Z_\alpha(x_1), Z_\alpha(x_2)) = \alpha_1(x_1, x_2) - \alpha(x_1)\alpha(x_2)$, $\mathrm{cov}(Z_\beta(x_1), Z_\beta(x_2)) = \beta_1(x_1, x_2) - \beta(x_1)\beta(x_2)$, for all $x_1, x_2 \in D$.*

LEMMA 2. *Let $x_1, \ldots, x_k \in D$. Under conditions $(A_{0,1})$ and $(A_4)$, for all $\boldsymbol{t} = (t_1, \ldots, t_k)' \in \mathbb{R}^k$, we have $\sum_{j=1}^k t_j \tilde{Z}_n(x_j) \xrightarrow{D} \mathrm{N}(0, \boldsymbol{t}'\Sigma\boldsymbol{t})$, where*

$$(\Sigma)_{jl} = \frac{\varphi_1(x_j, x_l)}{\psi(x_j)\psi(x_l)} + \frac{\varphi(x_j)\varphi(x_l)\psi_1(x_j, x_l)}{\psi^2(x_j)\psi^2(x_l)}$$

$$- \frac{\varphi(x_l)\mu(x_j, x_l)}{\psi(x_j)\psi^2(x_l)} - \frac{\varphi(x_j)\mu(x_j, x_l)}{\psi(x_l)\psi^2(x_j)}.$$

Put $Z_n(x) = \sqrt{n}(\hat{m}(x) - m(x)) = X_n(x) + Y_n(x)$, where $\psi(x)X_n = \sqrt{n}(\hat{\varphi}(x) - \varphi(x))$ and $\hat{\psi}(x)\psi(x)Y_n(x) = -\sqrt{n}(\hat{\psi}(x) - \psi(x))\hat{\varphi}(x)$. It suffices to prove (a) convergence of the finite-dimensional limit distribution of $Z_n$, and (b) tightness of $Z_n$. To establish (a), note that

(5.1) $$Z_n(x) = \tilde{Z}_n(x) - \sqrt{n}\frac{\hat{\psi}(x) - \psi(x)}{\psi(x)}\left[\frac{\hat{\varphi}(x)}{\hat{\psi}(x)} - \frac{\varphi(x)}{\psi(x)}\right].$$

Now

$$\sup_{x \in D} \left| \frac{\hat{\varphi}(x)}{\hat{\psi}(x)} - \frac{\varphi(x)}{\psi(x)} \right| \leq \frac{\sup_{x \in D} |\hat{\varphi}(x) - \varphi(x)|}{\inf_{x \in D} |\hat{\psi}(x)|}$$

$$+ \frac{\sup_{x \in D} |\varphi(x)| \cdot \sup_{x \in D} |\psi(x) - \hat{\psi}(x)|}{\inf_{x \in D} |\psi(x)\hat{\psi}(x)|},$$

where $\inf_{x \in D} |\psi_n(x)| \xrightarrow{P} \inf_{x \in D} |\psi(x)| > 0$, which, combined with Lemma 1, proves that the last term of (5.1) tends to zero as $n$ tends to infinity, and thus $Z_n(x)$ has the same finite-dimensional limit distribution as $\tilde{Z}_n(x)$. From Lemma 2 and the Cramér–Wold device, this limit distribution is the same as that claimed for $Z$ in Theorem 1. To prove (b), note that, by the proof of Lemma 1, the sequences $\sqrt{n}(\hat{\varphi}(x) - \varphi(x))$ and $\sqrt{n}(\hat{\psi}(x) - \psi(x))$ are tight. The sequence $\hat{\varphi}(x)/\hat{\psi}(x)$ is tight if we show that for given $\epsilon, \eta \geq 0$ and sufficiently small $\delta$ and large $n$,

$$(5.2) \qquad P\left( \sup_{|x-y| \leq \delta} |\hat{\varphi}(x)/\hat{\psi}(x) - \hat{\varphi}(y)/\hat{\psi}(y)| \geq \epsilon \right) \leq \eta.$$

Now, defining $\xi(x) = \iint y f_{W,Y}(w, y) f_\delta'(x - w) \, dw \, dy$, $\zeta(x) = \int f_W(w) f_\delta'(x - w) \, dw$, $\hat{\xi}(x) = n^{-1} \sum_{i=1}^n Y_i f_\delta'(x - W_i)$ and $\hat{\zeta}(x) = n^{-1} \sum_{i=1}^n f_\delta'(x - W_i)$, let $\hat{T}(x) = [\hat{\xi}(x)\hat{\psi}(x) - \hat{\varphi}(x)\hat{\zeta}(x)]/\hat{\psi}^2(x)$. By the mean value theorem, the left-hand side of (5.2) is bounded by $P(\sup_{x \in D} |\hat{T}(x)| \geq \epsilon/\delta)$ and (5.2) follows if we note that

$$\sup_{x \in D} |\hat{T}(x)| \leq \sup_{x \in D} |\hat{\xi}(x) - \xi(x)|/|\hat{\psi}(x)| + \sup_{x \in D} |\xi(x)|/|\hat{\psi}(x)|$$

$$+ \left[ \sup_{x \in D} |\hat{\varphi}(x) - \varphi(x)|/|\hat{\psi}(x)|^2 + \sup_{x \in D} |\varphi(x)|/|\hat{\psi}(x)|^2 \right]$$

$$\times \left[ \sup_{x \in D} |\hat{\zeta}(x) - \zeta(x)|/|\hat{\psi}(x)|^2 + \sup_{x \in D} |\zeta(x)|/|\hat{\psi}(x)|^2 \right],$$

which tends to zero as $n$ tends to infinity. Property (b) follows.

PROOF OF LEMMA 1. We prove the result for $\alpha$; the proof for $\beta$ is analogous. Let $x_1, \ldots, x_k \in D$, $\hat{\boldsymbol{\alpha}} = (\hat{\alpha}(x_1), \ldots, \hat{\alpha}(x_k))'$, $\boldsymbol{\alpha} = (\alpha(x_1), \ldots, \alpha(x_k))'$ and $\boldsymbol{Z}_\alpha \sim N_k(0, \Sigma_\alpha)$, where $(\Sigma_\alpha)_{ij} = \alpha_1(x_i, x_j) - \alpha(x_i)\alpha(x_j)$. Applying the central limit theorem to the i.i.d. sequence $T_1, \ldots, T_n$, with $T_i = \sum_{j=1}^k t_j Y_i h(x_j, W_i)$, it is not hard to prove that, for all $\boldsymbol{t} = (t_1, \ldots, t_k)' \in \mathbb{R}^k$, $\sqrt{n} \boldsymbol{t}'(\hat{\boldsymbol{\alpha}} - \boldsymbol{\alpha}) \xrightarrow{D} \boldsymbol{t}' \boldsymbol{Z}_\alpha$. From the Cramér–Wold device, we deduce that $\sqrt{n}(\hat{\boldsymbol{\alpha}} - \boldsymbol{\alpha}) \xrightarrow{D} \boldsymbol{Z}_\alpha$, which implies weak convergence of the finite-dimensional distributions. Using uniform Lipschitz continuity of $h$ in the first coordinate, one can show that $\mathrm{E}(\sqrt{n}[\hat{\alpha}(x_1) - \alpha(x_1) - \hat{\alpha}(x_2) + \alpha(x_2)])^2 \leq c|x_1 - x_2|^2$, which implies tightness of $\sqrt{n}(\hat{\alpha} - \alpha)$. □

PROOF OF LEMMA 2. Since $\mathrm{E}\hat{\varphi}(x) = \varphi(x)$ and $\mathrm{E}\hat{\psi}(x) = \psi(x)$, we have that $\mathrm{E}(\tilde{Z}_n(x)) = 0$. The result follows from the central limit theorem if we note that $\sqrt{n}\sum_{j=1}^{k} t_j \tilde{Z}_n(x_j)$ may be written as $\sum_{i=1}^{n} T_i$, where $T_i = \sum_{j=1}^{k} t_j f_\delta(x_j - W_i)[Y_i/\psi(x_j) - \varphi(x_j)/\psi^2(x_j)]$. $\square$

5.2. *Outline proof of Theorem* 2. We shall derive the second result at (2.4); a proof of the first result there is similar. Define the functions $a = (f_W g)^{\mathrm{ft}}$, $\hat{a} = \widehat{(f_W g)}^{\mathrm{ft}}$, $b = (f_\delta^{\mathrm{ft}})^2$, $\hat{b} = (\hat{f}_\delta^{\mathrm{ft}})^2$, $c = f_\delta^{\mathrm{ft}}$, $\Delta_a = \hat{a} - a$ and $\Delta_b = \hat{b} - b$. Let $\mathcal{T}$ denote the interval $[-\tau_n, \tau_n]$, write $\tilde{\mathcal{T}}$ for the complement in $\mathbb{R}$ of $\mathcal{T}$, and put $u_x(t) = e^{-itx}$. Then, uniformly in $x$,

$$
\begin{aligned}
2\pi\tilde{\varphi}(x) &= \int_{\mathcal{T}} \hat{a}|\hat{b}|^{1/2} u_x \\
&= \int_{\mathcal{T}} (a + \Delta_a)|b|^{1/2}(1 + b^{-1}\Delta_b)^{1/2} u_x \\
&= \int_{\mathcal{T}} (a + \Delta_a)cu_x + O_p\left[\int_{\mathcal{T}} |a/c|(\mathrm{E}\Delta_b^2)^{1/2} + \int_{\mathcal{T}} c^{-1}(\mathrm{E}\Delta_a^2 \mathrm{E}\Delta_b^2)^{1/2}\right].
\end{aligned}
$$

Using the fact that $\hat{a}$ equals a sum of $n$ independent and identically distributed random variables, and $\hat{b}$ is expressed in a form similar to a $U$-statistic, it can be shown that $\mathrm{E}[\Delta_a(t)^2] = O(n^{-1})$ and $\mathrm{E}[\Delta_b(t)^2] = O(N^{-1})$, uniformly in $t$. Moreover, (2.2) and (2.3) imply that $\int_{\mathcal{T}} |a/c| = O(1)$, $\int_{\mathcal{T}} c^{-1} = O(\tau_n^{\lambda_\delta + 1})$, $\int_{\tilde{\mathcal{T}}} acu_x = O(\tau_n^{1 - \lambda - \lambda_\delta})$ and $\int_{\tilde{\mathcal{T}}} \Delta_a cu_x = O_p(n^{-1/2}\tau_n^{1 - \lambda_\delta})$, the latter two results holding uniformly in $x$. Therefore, uniformly in $x$,

$$
\begin{aligned}
2\pi\tilde{\varphi}(x) &= \int (a + \Delta_a)cu_x \\
(5.3) \qquad &\quad + O_p(N^{-1/2} + n^{-1/2}N^{-1/2}\tau_n^{\lambda_\delta + 1} + \tau_n^{1 - \lambda - \lambda_\delta} + n^{-1/2}\tau_n^{1 - \lambda_\delta}) \\
&= \int \hat{a}cu_x + o_p(n^{-1/2}).
\end{aligned}
$$

Since $\hat{\varphi}(x) = (2\pi)^{-1}\int \hat{a}cu_x$, then the second part of (2.4) follows from (5.3).

5.3. *Case where* $X = p(T \mid \theta) + \delta$. This generalization, in which $(T, Y)$ rather than $(W, Y)$ is observed, was introduced in Section 1.3. There we noted that the unknown parameter $\theta$ could be estimated by least squares from data $(T_i', X_i')$, for $1 \le i \le r$, on $(T, X)$. In the case of a linear model, $p(t \mid \theta) = \theta^{(1)} + \theta^{(2)}t$, and our estimator of $W_i = \theta^{(1)} + \theta^{(2)}T_i$ is $\widehat{W}_i = \hat{\theta}^{(1)} + \hat{\theta}^{(2)}T_i$. We shall treat this particular case below; other models for $p$ can be addressed similarly.

Let $\hat{m}^*$, $\hat{\varphi}^*$, $\hat{\psi}^*$ and $\tilde{\psi}^*$ denote the versions of $\hat{m}$, $\hat{\varphi}$, $\hat{\psi}$ and $\tilde{\psi}$, respectively, obtained on replacing $W_i$ by $\widehat{W}_i$ throughout. It will be assumed that $n = O(r)$. In this case the least squares estimators $\hat{\theta}^{(1)}$ and $\hat{\theta}^{(2)}$ are $\sqrt{n}$-consistent.

First we consider the setting where $f_\delta$ is known. Provided $f_\delta$ has two bounded derivatives, we may write

$$
\hat{\varphi}^*(x) = n^{-1} \sum_{i=1}^n Y_i f_\delta(x - \widehat{W_i})
$$

(5.4)
$$
= \hat{\varphi}(x) - (\hat{\theta}^{(1)} - \theta^{(1)}) \mathrm{E}\{Y f_\delta'(x - W)\}
$$
$$
- (\hat{\theta}^{(2)} - \theta^{(2)}) \mathrm{E}\{TY f_\delta'(x - W)\} + o_p(n^{-1/2}),
$$

$$
\hat{\psi}^*(x) = n^{-1} \sum_{i=1}^n f_\delta(x - \widehat{W_i})
$$

(5.5)
$$
= \hat{\psi}(x) - (\hat{\theta}^{(1)} - \theta^{(1)}) \mathrm{E}\{f_\delta'(x - W)\}
$$
$$
- (\hat{\theta}^{(2)} - \theta^{(2)}) \mathrm{E}\{T f_\delta'(x - W)\} + o_p(n^{-1/2}).
$$

Here, $\hat{\varphi}$ and $\hat{\psi}$ are the original estimators of $\varphi$ and $\psi$ given in Section 1.2 for the case where $W_i$ is directly observed; $W = \theta^{(1)} + \theta^{(2)} T$; and the remainder terms $o_p(n^{-1/2})$ are uniform in $x$, provided the conditions of Theorem 2 hold and, in addition, $\mathrm{E}(T^2) < \infty$.

It follows from (5.4) and (5.5) that $\hat{\varphi}$ and $\hat{\psi}$ are $\sqrt{n}$-consistent for $\phi$ and $\psi$, respectively, and $\hat{m}^* = \hat{\varphi}^*/\hat{\psi}^*$ is $\sqrt{n}$-consistent for $m$. A version of Theorem 1 is readily obtained in this setting, using (5.4) and (5.5). Unless $r/n \to \infty$, the covariance structure of the limiting Gaussian process depends on whether the data $(T_i', X_i')$, from which $\hat{\theta}^{(1)}$ and $\theta^{(2)}$ are computed, are independent of the data $(W_i, Y_i)$ used to calculate $\hat{\varphi}$ and $\hat{\psi}$, or whether $(T_i', X_i') = (T_i, X_i)$ and the triples $(T_i, X_i, Y_i)$ are observed.

The case where $f_\delta$ is not known, and is consistently estimated from replicated data as discussed in Section 1.2, is similar although more complex. Our estimator $\hat{f}_\delta^{\mathrm{ft}}$, given at (1.6), does not alter since it does not use the data $W_i$. On the other hand, the estimators $\hat{f}_W^{\mathrm{ft}}$ and $\widehat{(f_W g)}^{\mathrm{ft}}$, given at (1.6) and (1.7), are replaced by

$$
\hat{f}_W^{\mathrm{ft}\,*}(t) = n^{-1} \sum_{j=1}^n \exp(it\widehat{W_j}),
$$

$$
\widehat{(f_W g)}^{\mathrm{ft}\,*} = n^{-1} \sum_{j=1}^n Y_j \exp(it\widehat{W_j}).
$$

Substituting the latter for $\hat{f}_W^{\mathrm{ft}}$ and $\widehat{(f_W g)}^{\mathrm{ft}}$, respectively, in (1.8); Taylor-expanding $\exp(-it\widehat{W_j})$ as $\exp(itW_j)\{1 + it(\widehat{W_j} - W_j) + \cdots\}$; and taking the smoothing parameter $\tau_n$ in (1.8) to be of order $n^{(1/2)-2\eta}$, for some $\eta > 0$ [so that, under moment conditions on $W_j$, $\tau_n \sup_{j \le n} |\widehat{W_j} - W_j| = O_p(n^{-\eta})$], we may deduce that (2.4) continues to hold if $\hat{\varphi}$ and $\hat{\psi}$ there are replaced by $\hat{\varphi}^*$ and $\hat{\psi}^*$, provided more restrictive assumptions than those given in Theorem 2 are imposed.

**Acknowledgment.** We are grateful to Chris Klaassen for helpful discussion, and to the Editor, the Associate Editor and two reviewers for valuable comments and suggestions which led to considerable improvement of the paper.

A. DELAIGLE
DEPARTMENT OF MATHEMATICS
UNIVERSITY OF BRISTOL
BRISTOL BS8 1TW
UNITED KINGDOM
E-MAIL: aurore.delaigle@bristol.ac.uk

P. HALL
DEPARTMENT OF MATHEMATICS
 AND STATISTICS
UNIVERSITY OF MELBOURNE
PARKVILLE, VICTORIA 3010
AUSTRALIA
E-MAIL: hall@unimelb.edu.au

H.-G. MÜLLER
DEPARTMENT OF STATISTICS
UNIVERSITY OF CALIFORNIA
DAVIS, CALIFORNIA 95616
USA
E-MAIL: mueller@wald.ucdavis.edu